\documentclass[12pt,reqno]{amsart}
\usepackage{amssymb}
\usepackage{cases}
\usepackage{bbm}
\usepackage{mathabx}
\usepackage{mathrsfs}
\usepackage{graphicx}
\usepackage{amsfonts}
\usepackage{enumerate}

\usepackage{amssymb, amsmath, amsfonts, latexsym}
\usepackage{enumerate}

\newtheorem{thm}{Theorem}[section]
\newtheorem{cor}[thm]{Corollary}
\newtheorem{lem}[thm]{Lemma}

\newtheorem{remarks}[thm]{Remark}

\theoremstyle{definition}
\newtheorem{definition}{Definition}[section]
 \theoremstyle{remark}

\usepackage{color}

\newcommand{\ee}{\mathbb{E}}

\newcommand{\rr}{\mathbb{R}}
\newcommand{\pp}{\mathbb{P}}
\newcommand{\qq}{\mathbb{Q}}

\def\BB{\mathcal B}
\def\CC{\mathcal C}

\def\FF{\mathcal F}
\def\EE{\mathbb E}

\def\SS{\mathcal S}

\def\BB{\mathcal B}
\def\CC{\mathcal C}

\def\FF{\mathcal F}
\def\EE{\mathbb E}
\def\HH{\mathcal H}

\def\SS{\mathcal S}

\def\<{\langle}
\def\>{\rangle}

\def\beq{\begin{equation}}
\def\nneq{\end{equation}}

\def\bdef{\begin{defn}}
\def\ndef{\end{defn}}

\def\bthm{\begin{thm}}
\def\nthm{\end{thm}}

\def\bprop{\begin{prop}}
\def\nprop{\end{prop}}

\def\brmk{\begin{remarks}}
\def\nrmk{\end{remarks}}

\def\bexa{\begin{exa}}
\def\nexa{\end{exa}}

\def\blem{\begin{lem}}
\def\nlem{\end{lem}}

\def\bcor{\begin{cor}}
\def\ncor{\end{cor}}

\def\<{\langle}
\def\>{\rangle}

\date{}

\def\bexe{\begin{exe}}
\def\nexe{\end{exe}}

\def\bprf{\begin{proof}}
\def\nprf{\end{proof}}

\def\bdes{\begin{description}}
\def\ndes{\end{description}}

\title[Transportation cost-information inequalities for stochastic wave equation ]{Transportation cost-information inequality for stochastic wave equation}

\author{Yumeng Li}
\address{Yumeng Li \\School of Statistics and Mathematics, Zhongnan University of Economics and Law, 430073, P. R. China.}
\email{li-yu-meng@163.com}

\author{XinYu Wang}
\address{XinYu Wang \\School of  Mathematics and Statistics, Huazhong University of Science and Technology, 430073, P. R. China.}
\email{wang\_xin\_yu@hust.edu.cn}
\date{}

\begin{document}
\maketitle

 \noindent {\bf Abstract:}
 In this paper, we prove a Talagrand's $T_2$ transportation cost-information inequality  for the law of  a  stochastic wave equation in spatial dimension  $d=3$ driven by the Gaussian random field, white in time and correlated in space, on the continuous paths space with respect to the uniform topology.

 \vskip0.3cm

 \noindent{\bf Keyword:} {Stochastic wave equation; Girsanov transformation; Transportation cost-information inequality.
}
 \vskip0.3cm

\noindent {\bf MSC: } {60H15; 60H20.}
\vskip0.3cm

\section{Introduction}

    The purpose of this paper is to study the Talagrand's $T_2$ transportation cost-information inequality for the   following   stochastic wave equation in spatial dimension  $d=3$:
\begin{equation}\label{SPDE}
    \begin{cases}
     \left(\frac{\partial^2}{\partial t^2}-\Delta\right) u(t,x)=  b\big(u(t,x)\big)+\dot{F}(t,x),\\
      u(0,x)=\nu_{1}(x),\\
     \frac{\partial}{\partial t}u(0,x)=\nu_{2}(x)
    \end{cases}
\end{equation}
for all $(t,x)\in [0,T]\times\rr^3$,
where  the coefficient $ b:\rr \rightarrow\rr $ is Lipschitz continuous, the term $\Delta u$ denotes the Laplacian of $u$ in the $x$-variable and the process $\dot{F}$ is the formal derivative of a Gaussian random field, white in time and correlated in space.
We recall that a random field solution to \eqref{SPDE} is a family of random variables $\{u(t,x), t\in\mathbb R_+, x\in \mathbb R^3\}$ such that $(t,x)\mapsto u(t,x)$ from $\mathbb R_+\times \mathbb R^3$ into $L^2(\Omega)$ is continuous and solves an integral form of \eqref{SPDE}, see Section \ref{sec framework} for details.

It is known that random field solutions have been shown to exist when $d\in\{1,2,3\}$, see \cite{Dalang1999}. In spatial dimension $1$,  a solution to the non-linear wave equation driven by space-time white noise was given in \cite{Walsh} by using Walsh's martingale measure stochastic integral. In dimensions $2$ or higher, there is no function-valued solution with space-time white noise, some spatial correlation is needed.  A necessary and sufficient condition on the spatial correlation for existence of a random field solution was given in \cite{Dalang Frangos}.
 Since the fundamental solution in spatial dimension $d=3$ is not a function, this required an extension of Walsh's martingale measure stochastic integral to integrands that are Schwartz distributions, the existence of a random field solution to \eqref{SPDE} is given in \cite{Dalang1999}.  H\"older continuity of the solution was   established in  \cite{DS2009}. In spatial dimensional $d\ge4$, since the fundamental solution of the wave equation is not a measure, but a Schwarz distribution  that is a derivative of some order of a measure, the methods used in dimension $3$ do not apply to higher
dimensions, see \cite{Dalang Frangos} for the study of the solutions.

  Transportation cost-information inequalities have been recently deeply studied, especially
for their connection with the concentration of measure phenomenon,  log-Sobolev inequality, Poincar\'e inequalities and  Hamilton-Jacobi's equation, see \cite{BGL2014, BGL2001, BLM, Goz, Ledoux, OV,Tal, Villani}.

Let us   recall the transportation inequality. Let $(E,d)$ be a metric space equipped with $\sigma$-field $\mathcal B$ such that $d(\cdot, \cdot)$ is $\mathcal B\times \mathcal B$ measurable. Given $p\ge1$ and two probability measures $\mu$ and $\nu$ on $E$, the Wasserstein distance is defined by
 $$
 W_{p, d}(\mu, \nu):=\inf_{\pi}\left[\int d^p(x,y)\pi(dx,dy) \right]^{\frac1p},
 $$
where the infimum is taken over all  the probability measures  $\pi$   on $E\times E$ with marginal distributions $\mu$ and $\nu$.  The relative entropy of $\nu$ with respect to (w.r.t. for short) $\mu$ is defined as
\begin{equation}\label{entropy}
\mathbf H(\nu|\mu):=\left\{
       \begin{array}{ll}
         \int \log \frac{d\nu}{d\mu}d\nu,   & \text{if } \nu\ll \mu ;\\
        +\infty, & \hbox{\text{otherwise}.}
       \end{array}
     \right.
\end{equation}
\begin{definition}
The probability measure $\mu$ is said to  satisfy the   transportation cost-information inequality ${\bf T_{p}}(C)$ on $(E, d)$ if there exists a constant $C>0$ such that for any probability measure $\nu$ on $E$,
$$
W_{p,d}(\mu,\nu)\le \sqrt{2C {\mathbf H}(\nu|\mu)}.
$$
\end{definition}
  Recently, the problem of  transportation inequalities and their applications to diffusion  processes have been widely studied.  The ${\bf T_2}(C)$ inequality, first established by M. Talagrand \cite{Tal} for the Gaussian measure with the sharp constant $
C=2$.    The approach of M. Talagrand is generalized  by D. Feyel and
A.S. \"Ust\"unel \cite{FU}  on the abstract Wiener space with respect to Cameron-Martin distance  using the Girsanov theorem.
   With regard to the paths of finite   stochastic differential equation (SDE for short), by means of Girsanov transformation and the martingale representation theorem, the ${\bf T_2}(C)$  w.r.t. the $L^2$ and the Cameron-Martin distances were established by  H.  Djellout {\it et al.} \cite{DGW}; the ${\bf T_2}(C)$ w.r.t. the uniform metric was obtained by \cite{Ust,WZ2004}.  J. Bao {\it et al.} \cite{BWY} established the ${\bf T_2}(C)$ w.r.t. both the uniform and the $L^2$ distances on the path space for the segment process associated to a class of neutral function stochastic differential equations.  B. Saussereau \cite{Sau} studied the ${\bf T_2}(C)$ for SDE driven by a fractional Brownian motion, and  S. Riedel \cite{Rie} extended this result to the law of SDE driven by general Gaussian processes by using Lyons' rough paths theory.
S. Pal \cite{Pal} proved that probability laws of certain multidimensional semimartingales
which includes time-inhomogenous diffusions,  satisfy
quadratic transportation cost inequality under the uniform metric.  Those,  in particular,  imply some results about concentration of boundary local time of reflected Brownian motions.

Motivated  by the source of the noise modeled by the random terms in  partial differential equations, which include physical noise (such as thermal noise), the stochastic partial differential equations have been studied in many literatures in past thirty years. For the  stochastic reaction-diffusion equation, L. Wu and Z. Zhang \cite{WZ2006} studied the ${\bf T_2}(C)$ w.r.t. $L^2$-metric by Galerkin's approximation. By Girsanov's transformation, B. Boufoussi and S. Hajji \cite{BH} obtained the ${\bf T_2}(C)$ w.r.t. $L^2$-metric for the stochastic heat equations driven by space-time white noise and driven by fractional noise.  Those  results are established for the stochastic parabolic equations.  However,  the hyperbolic case is   much more complicated, one difficulty comes from the  more complicated stochastic integral, another one comes from the lack of good regularity properties of  the  fundamental solutions. See  \cite{Dalang}  for the study of the stochastic wave equation.

In this paper, we shall study the  Talagrand's $T_2$-transportation inequality  for the law of  a  stochastic wave equation \eqref{SPDE} on the continuous paths space with respect to the  uniform metric.

The rest of this paper is organized as follows. In Section 2, we first give  the properties of the  Eq. \eqref{SPDE}, and then state the main  result of this paper.  In Section 3, we shall prove the main result.

\section{Framework and main result}\label{sec framework}
\subsection{Framework}
For any
$d\geq 1$, let $\SS(\rr^{d+1})$ be the space of Schwartz functions, all of whose derivatives are rapidly decreasing.  $F=(F(\varphi),\varphi \in \SS(\rr^{d+1}))$  is a  Gaussian process  defined on some probability space with zero mean and covariance functional
\begin{equation}\label{covariance}
\ee(F(\varphi)F(\psi))=J(\varphi, \psi):=\int_{\rr_+}ds\int_{\rr^d}dx\int_{\rr^d}dy  \varphi(s, x)f(x-y){\psi}(s,y),
\end{equation}
where $f:\mathbb R^d\rightarrow\rr_+$  continuous on $\mathbb R^d\setminus\{0\}$.

 According to \cite{Dalang1999}, there are some requirements on $f$.  As a covariance functional of a Gaussian process, the function $J(\cdot, \cdot)$ should be non-negative definite, this implies that $f$ is symmetric ($f(x)=f(-x)$ for all $x\in \mathbb R^d$), and is equivalent to the existence of a non-negative tempered measure $\mu$ on $\mathbb R^d$, whose Fourier transform is $f$. More precisely, for any $\varphi\in \mathcal S(\mathbb R^d)$, let $\mathcal F\varphi$ be the Fourier transform of $\varphi$:
 $$
 \mathcal F\varphi(\xi):=\int_{\mathbb R^d}\exp(-2i\pi\xi\cdot x)\varphi(x)dx.
 $$
The relationship between $\mu$ and $f$ is, by definition of the Fourier transform on the space $\mathcal S'(\mathbb R^d)$ of tempered distributions, that is for all $\varphi\in \mathcal S(\mathbb R^d)$,
$$
\int_{\mathbb R^d} f(x)\varphi(x)dx=\int_{\mathbb R^d}\mathcal F\varphi(\xi)\mu(d\xi).
$$
Elementary properties of Fourier transform show that for all $\varphi, \psi\in \mathcal S(\mathbb R^d)$,
\begin{equation}\label{eq HH}
\begin{split}
 \langle\varphi,\psi \rangle_\HH&:=\int_{\mathbb R^d}dx\int_{\rr^d}dy  \varphi(x)f(x-y){\psi}(y) \\
&= \int_{\mathbb R^d}\mu(d\xi)\mathcal F\varphi(\xi)\overline{\mathcal F\psi(\xi)}.
\end{split}
\end{equation}
Here $\bar z$ is the complex conjugate of $z$.

According to \cite{DM}, the Gaussian   process $F$ with covariance \eqref{covariance} can be extended to a   martingale measure
$$
M=\left\{M_t(A),\ t\ge0,\ A\in\BB_b(\rr^3)\right\},
$$
where $\BB_b(\rr^3)$ denotes the collection of all bounded Borel measurable sets in $\rr^3$.

 For each $t\ge0$, denote by $\mathcal F_t$ the $\sigma$-field generated by the random variables $\{M_s(A), s\in[0,t],A\in \mathcal B(\rr^3)\}$, that is
 $$
 \mathcal F_t:=\sigma\{M_s(A), s\in[0,t],A\in \mathcal B(\rr^3)\}.
 $$

Let $\HH$ be the Hilbert space obtained by the completion of $\SS(\rr^3)$ with the inner product $\langle\cdot, \cdot\rangle_\HH $ defined by \eqref{eq HH}, and denote by $\|\cdot\|_{\HH}$ the induced norm. Let $\HH_T:=L^2([0,T]; \HH)$ and consider the usual $L^2$-norm  $\|\cdot\|_{\HH_T}$ on this space. Then $\HH_T$ is a Hilbert space with the inner product
\beq\label{eq HHT}
\langle \psi_1,\psi_2 \rangle_{\HH_T}:=\int_0^T\langle \psi_1(t),\psi_2(t) \rangle_\HH dt, \ \ \ \varphi,\psi\in \HH_T.
\nneq

Let $D$ be an open set $D$ in Euclidean space $\mathbb R^d$ for $d\ge1$.  For any integer  $n\ge1$,  let $\CC^n(D)$ be the space of all continuous  functions from $D$ to $\rr$, whose derivatives up to  order $n$  are also continuous. For any $\delta\in(0,1)$,  let $\CC^{\delta}(D)$ be the space of all H\"{o}lder continuous functions of degree $\delta$, with the H\"older norm
$$
\|g\|_{\delta}:=\sup_{x\neq y}\frac{|g(x)-g(y)|}{|x-y|^\delta}<\infty, \ \forall g \in \CC^{\delta}(D);
$$
 and let $\CC_{Lip}(D)$ be the space of all Lipschitz continuous functions, with the norm
$$
\|f\|_{Lip}:=\sup_{x\neq y}\frac{|g(x)-g(y)|}{|x-y|}<\infty,\ \forall g \in \CC_{Lip}(D).
$$
\vskip0.3cm

{\bf Hypothesis (H)}:
\begin{itemize}
        \item[(H.1)]
      There exists a constant $K>0$ such that
\beq\label{Lip}
 |b(x)-b(y)|\le K|x-y|, \quad\forall x,y\in\rr.
\nneq

        \item[(H.2)] The function $f$ given in \eqref{covariance} can be expressed by            $f(x)=\varphi(x)|x|^{-\beta},x\in \rr^3 \backslash \{0\}$, with $\beta\in]0,2[$. Here the functions $\varphi$ and $\nabla \varphi$ are bounded,
            $0<\varphi\in\CC^1(\rr^3), \nabla\varphi\in\CC^\delta(\rr^3)$ with  $\delta\in]0,1]$.
         \item[(H.3)] The initial values $\nu_1,\nu_2$ are bounded, $\nu_1\in \CC^2(\rr^3),\nabla \nu_1 $ is bounded,  $\triangle\nu_1$ and $\nu_2$ are H\"{o}lder continuous with degrees $\gamma_1,\gamma_2\in]0,1]$, respectively.
\end{itemize}
\vskip0.3cm

We remark that the  hypothesis (H.2) implies that, for any $T>0$,
\beq\label{eq G}
M(T):= \sup_{t\in[0,T]}\int_{\rr^3}|\FF G(t)(\xi)|^2\mu(d\xi)<\infty,
 \nneq
 see \cite{Dalang}.
\vskip0.3cm

By Walsh's theory of stochastic integration with respect to (w.r.t. for short)  martingale measures, for any $t\geq 0$ and  $h \in \HH$, the stochastic integral
$$
B_t(h):=\int_0^t\int_{\rr^d}h(y)M(ds,dy)
$$
is well defined, and
$$
\left\{B_t^k:=\int_0^t\int_{\rr^d}e_k(y)M(ds,dy);\ k\ge1\right\}
$$
defines a sequence of independent standard Wiener processes, here $\{e_k\}_{k\ge1}$ is a complete orthonormal system of the Hilbert space $\HH$.  Thus, $B_t:=\sum_{k\ge1}B_t^k e_k$ is a cylindrical Wiener process on $\HH$. See \cite{Walsh}.

According to Dalang and Sanz-Sol\'e \cite{DS2009}, under hypothesis ({\bf H}), Eq. \eqref{SPDE} admits a unique solution $u$:
\begin{equation}\label{SPDE solution}
\begin{split}
u(t,x)&=w(t,x)+ \sum_{k\ge 1}\int_0^t\left\langle G(t-s,x-\cdot), e_k(\cdot)\right\rangle_{\HH}
dB_s^k\\
           &+\int_0^t\int_{\mathbb R^3}\big[G(t-s,x,y) b(u(s,y)) dyds,
           \end{split}
\end{equation}
 where
$$
w(t,x):= \frac{d}{dt}G(t,x,y)\nu_1(y)  + G(t,x,y)\nu_2(y),
$$
with $G(t,x,y)=\frac{1}{4\pi t}\sigma_t(x-y) $,   $\sigma_t$ is the uniform surface measure (with total mass $4\pi t^2$) on the sphere of radius $t$.
 Furthermore, for any $p\in[2,\infty[$,
\beq\label{eq u e estimate}
\sup_{(t,x)\in[0,T]\times\rr^3}\ee\left[|u(t,x)|^p\right]<+\infty.
\nneq
   See Dalang and Sanz-Sol\'e \cite{DS2009}  or  Hu {\it et al} \cite{HHN} for details.

\subsection{Main results}

Let $\mathcal C([0,T]\times \rr^3)$ be the space of all continuous functions from $[0,T]\times \mathbb R^3$ to $\mathbb R$, endowed with the uniform norm
$$
\|f\|_{\infty}:=\sup_{(t,x)\in [0,T]\times\mathbb R^3}|f(t,x)|.
$$

For initial function $\nu:=(\nu_1, \nu_2)$ satisfied  (H.3), let $\pp_{\nu}$ be the law of $\{u(t,x), (t,x)\in[0,T]\times \rr^3\}$ on $\mathcal C([0,T]\times \rr^3)$ with initial value $u(0,x)=\nu_1(x)$ and   $\frac{\partial}{\partial t}u(0,x)=\nu_{2}(x)$.

 Recall the constants $K, M(T)$ given by \eqref{Lip} and \eqref{eq G} respectively.  In this paper, we establish the following result:
 \bthm\label{thm transport}  Under  Hypothesis (H), there exists a constant $C(T, K):= TM(T)e^{\frac{T^4K^2}{2}}$ such that the probability measure $\pp_{\nu}$ satisfies ${\bf T_2}(C)$ on the space $\mathcal C([0,T]\times \rr^3)$ endowed with the uniform norm $\|\cdot\|_{\infty}$.

 \nthm

As indicated in \cite{BGL2001}, many interesting consequences can be derived from Theorem \ref{thm transport},
see also Corollary 5.11 of \cite{DGW}.
\bcor
Under  Hypothesis (H),   we have for any $T>0$, the following statements hold for the constant $C(T, K)=TM(T)e^{\frac{T^4K^2}{2}}$.
\begin{itemize}
\item[(a)] For any Lipschitzian function $U$ on $\mathcal C([0,T]\times \rr^3)$ with respect to uniform norm, we have
\begin{equation*}
 \mathbb E^{\mathbb P_{\nu}}\exp\left(U-\ee^{\pp_{\nu}}U\right)\le e^{\frac {C(T, K)}2 \|U\|_{\text {Lip}}^2}.
     \end{equation*}

  \item[(b)] (Inequality of Hoeffding type) For any $V:\rr\rightarrow\rr$ such that $\|V\|_{Lip} <\infty$, we have that for any $r\ge0$,
     \begin{equation*}
     \begin{split}
      \pp\left(\frac{1}{T}\int_0^T  \|V(u(t))\|_{\infty}dt-\ee\left[ \frac{1}{T}\int_0^T \|V(u(t))\|_{\infty}dt\right]>r \right)
      \le \exp\left(-\frac{r^2}{2C(T, K) \|V\|_{\text {Lip}}^2}\right).
      \end{split}
     \end{equation*}

  \end{itemize}

\ncor

\section{The proof}
\subsection{An important lemma}

 This lemma is an analogue of the result \cite[Theorem
5.6]{DGW} for the space-colored time-white noise instead of the finite-dimensional Brownian motion.

\blem\cite[Lemma 6.2]{KS}\label{Girsanov}
For every probability measure $\mathbb Q\ll \mathbb P_{\nu}$ on the space  $L^2([0,T]\times \mathbb R^3;\mathbb R)$, there exists an adapted $\mathbb Q$-a.s. $h=\{h(s,x), (s,x)\in[0,T]\times \rr^3\}$ such that $\|h\|_{\HH_T}<\infty$, $\mathbb Q$-a.s., and  the function $\widetilde F: L^2([0,T]\times \mathbb R^3;\mathbb R)\rightarrow L^2(\Omega)$ defined by
\beq
\tilde F(\phi):=F(\phi)-\int_0^{t}\int_{\mathbb R^3}  \phi(s,x) h(s,x)dx ds,
\nneq
is a space-colored time-white noise with the  spectral density $f$ with respect to the measure $\qq$.
The Randon-Nikodym derivative is given by
\beq\label{eq RN}
\frac{d\mathbb Q}{d\mathbb P_{\nu}}\Big|_{\mathcal F_T}=\exp\left(\int_0^T\int_{\rr^3} h(s,x)F(ds,dx)-\frac12\int_0^T\| h(s)\|^2_{\HH}ds \right),
\nneq
and the relative entropy is given by
\beq\label{eq entropy}
\mathbf H(\mathbb Q| \mathbb P_{\nu})=\frac12 \ee^{\mathbb Q}\left[\|h\|_{\HH_T}^2\right].
\nneq

\nlem

\subsection{The proof of Theorem \ref{thm transport}}
\bprf
  It is enough to prove the result for any probability measure $\qq$ on $\mathcal C([0,T]\times \rr^3)$
 such that $\qq\ll \pp_{\nu}$ and $\mathbf H(\qq|\pp_{\nu})<\infty$.  Assume that the Randon-Nikodym derivative is given by
\beq
\frac{d\mathbb Q}{d\mathbb P_{\nu}}\Big|_{\mathcal F_T}=\exp\left(\int_0^T\int_{\rr^3} h(s,x)F(ds,dx)-\frac12\int_0^T\| h(s)\|^2_{\HH}ds \right),
\nneq
and the relative entropy is
\beq
\mathbf H(\mathbb Q| \mathbb P_{\nu})=\frac12 \ee^{\mathbb Q}\left[\|h\|_{\HH_T}^2\right].
\nneq

  Let $(\Omega, \mathcal F, \widetilde \pp)$ be a complete probability space on which  $F$ is space-colored time-white noise with the  spectral density $f$ with respect to the measure $\qq$. Let
 $$\mathcal F_t=\mathcal F_t^{F}=\sigma(F(s, A), s\le t, \ \forall A\subset \mathbb R ^3)^{\widetilde \pp} \ \ \ \text{completion by } \widetilde \pp.$$
 Let $u_t(\nu)$ be the unique solution of \eqref{SPDE} with initial condition $\nu=(\nu_1, \nu_2)$. Then the law of $u_{\cdot}(\nu)$ is $\mathbb P_{\nu}$. Consider
 $$
 \widetilde {\mathbb Q}:=\frac{d\mathbb Q}{d\mathbb P_{\nu}} (u_\cdot(\nu))\cdot\widetilde {\mathbb P}.
 $$
  Then
  $$
\mathbf H(\mathbb Q|\mathbb P_{\nu})=\mathbf H(\widetilde{\mathbb Q}|\widetilde{\mathbb P})=\frac12 \ee^{\mathbb Q}\left[\|h\|_{\HH_T}^2\right]
$$

For a complete orthonormal system $\{e_k\}_{k\ge1}$ of the Hilbert space $\HH$,  let
$$
\left\{B_t^k:=\int_0^t\int_{\rr^3}e_k(y)F(ds,dy);\ k\ge1\right\}.
$$
Then $B_t:=\sum_{k\ge1}B_t^k e_k$ is a cylindrical Wiener process on $\HH$ under $\widetilde{\mathbb P}$, and $\sum_{k\ge1}(B_t^k+\langle h, e_k\rangle_{\HH}) e_k$ is a cylindrical Wiener process on $\HH$ under $\widetilde{\mathbb Q}$ by Lemma \ref{Girsanov}.

 According to Lemma \ref{Girsanov},   we couple $(\pp, \qq)$ as the law of a process $(u,v)$ under $\qq$
\begin{equation}\label{eq u}
\begin{split}
 u(t,x)&=w(t,x)+\int_0^t\int_{\rr^3} G(t-s,x,y) \widetilde{F}(ds,dy) \\
 &+\int_0^t\int_{\rr^3} G(t-s,x,y)b(u(s,y))dyds \\
 &+\sum_{k\ge 1}\int_0^t\left\langle G(t-s,x-\cdot), e_k(\cdot)\right\rangle_{\HH}\cdot\left\langle h(s,\cdot), e_k(\cdot)\right\rangle_{\HH}
ds,
\end{split}
                  \end{equation}
and
\begin{equation}\label{eq v}
\begin{split}
 v(t,x)&=w(t,x)+\int_0^t\int_{\rr^3} G(t-s,x,y) \widetilde F(ds,dy) \\
 &+\int_0^t\int_{\rr^3} G(t-s,x,y)b(v(s,y))dyds.
 \end{split}
                  \end{equation}
By the definition of the Wasserstein distance,
\beq\label{eq uv}
W_{2, \|\cdot\|_{\infty}}^2(\qq, \pp) \le \EE^{\widetilde\qq}\left[\sup_{(t,x)\in [0,T]\times \mathbb R^3} |u(t,x)-v(t,x)|^2 \right].
\nneq

In view of \eqref{eq entropy} and \eqref{eq uv}, it remains to prove that
\beq\label{eq tar}
\EE^{\widetilde\qq}\left[\sup_{(t,x)\in [0,T]\times \mathbb R^3} |u(t,x)-v(t,x)|^2 \right]\le C \EE^{\widetilde\qq}[\|h\|^2_{\HH_T}].
\nneq

From \eqref{SPDE solution},  \eqref{eq u} and \eqref{eq v}, we can represent $u(t,x)-v(t,x)$ as
\begin{equation}\label{T}
\begin{split}
 u(t,x)-v(t,x)&= \int_0^t\int_{\rr^3} G(t-s,x,y)[b(u(s,y))-b(v(s,y))]dyds \\
 &+\sum_{k\ge 1}\int_0^t\left\langle G(t-s,x-\cdot), e_k(\cdot)\right\rangle_{\HH}\cdot \left\langle  h(s,\cdot), e_k(\cdot)\right\rangle_{\HH}
ds \\
 =:& I_1(t,x)+I_2(t,x).
 \end{split}
   \end{equation}
   For every $t\in[0,T]$, define
$$
\eta(t)=\sup_{(s,x)\in[0,t]\times\rr^3} |u(s,x)-v(s,x)|^2.
$$
    By using the elementary inequality $(a+b)^2\le 2(a^2+b^2)$, we have
 \begin{align}\label{eq T}
 |u(t,x)-v(t,x)|^2\le &  2I_1^2(t,x)+2I_2^2(t,x).
   \end{align}
By the Cauchy-Schwarz inequality with with respect to the finite measure $G(t-s,x,y)dyds$ on $[0,T]\times \mathbb R^3$, with total measure $t^2/2$, and by the Lipschitz continuity of $b$, we obtain that for any $t\le T$,
\begin{equation}\label{T1}
\begin{split}
 |I_1(t,x)|^2 &\le K^2  \int_0^t\int_{\rr^3} G(t-s,x,y) dyds \cdot\int_0^t\int_{\rr^3} G(t-s,x,y) |u(s,y)- v(s,y) |^2dyds  \\
 &\leq     \frac{t^2 K^2}{2}\int_0^t\int_{\rr^3} G(t-s,x,y) \eta(s)dyds  \\
 &=  \frac{t^2 K^2}{2}\int_0^t(t-s) \eta(s)ds\\
 &\le \frac{T^3K^2}{2}\int_0^t \eta(s)ds,
 \end{split}
  \end{equation}
  where $\int_{\rr^3} G(t-s,x,y) dy=t-s$ is used in the last second line.  Let us estimate  the second term. By the Cauchy-Schwarz inequality and  \eqref{eq G}, we have for any $t\le T$,
\begin{equation}\label{T2}
\begin{split}
 |I_2(t,x)|^2
 &\leq   \int_0^t \|G(t-s, x, \cdot)\|_{\HH}^2ds \times  \int_0^t  \|h(s)\|_{\HH}^2ds
   \\
       &\leq  TM(T) \int_0^t  \|h(s)\|_{\HH}^2ds.
       \end{split}
\end{equation}

Putting \eqref{eq T}, \eqref{T1} and \eqref{T2} together, we have for any $t\le T$,
\begin{align*}
\eta(t)\le \frac{T^3K^2}{2}\int_0^t\eta(s)ds+TM(T) \|h\|_{\HH_T}^2.
\end{align*}
 Using the Gronwall's inequality, we obtain that
$$
 \ee^{\widetilde \qq}\left[ \eta(T) \right]\le  TM(T)e^{\frac{T^4K^2}{2}} \ee^{\widetilde \qq}[\|h\|_{\HH_T}^2].
 $$
The proof is complete.
  \nprf

{\bf Acknowledgements.} The authors are grateful to the anonymous referees for  comments and corrections.

\end{document}